\documentclass[a4paper,10pt]{amsart}

\usepackage[top=2.5cm, bottom=2.5cm, left=2.5cm, right=2.5cm]{geometry}

\usepackage{amssymb}
\usepackage{ascmac}
\usepackage[dvipdfmx]{graphicx}
\usepackage{color}
\usepackage{bbm}

\usepackage[dvipdfm,colorlinks=true,bookmarks=true,
bookmarksnumbered=true,bookmarkstype=toc,linktocpage=true
]{}

\newtheorem{Def}{Definition}[section]
\newtheorem{Thm}[Def]{Theorem}

\newtheorem{Assumption}[Def]{Assumption}
\newtheorem{Rem}[Def]{Remark}

\numberwithin{equation}{section}

\newcommand{\argmin}{\operatornamewithlimits {argmin}}

\allowdisplaybreaks

\usepackage{ascmac}

\newcommand{\mca}{\mathcal{A}}

\newcommand{\mcf}{\mathcal{F}}
\newcommand{\mci}{\mathcal{I}}
\newcommand{\mcl}{\mathcal{L}}
\newcommand{\mcm}{\mathcal{M}}

\newcommand{\mbbg}{\mathbb{G}}

\newcommand{\mbbn}{\mathbb{N}}

\newcommand{\mbbr}{\mathbb{R}}

\newcommand{\mbX}{\mathbf{X}}
\newcommand{\al}{\alpha}
\newcommand{\del}{\delta}
\newcommand{\ep}{\epsilon}

\newcommand{\D}{\Delta}
\newcommand{\sig}{\sigma}
\newcommand{\Sig}{\Sigma}

\newcommand{\gam}{\gamma}

\newcommand{\p}{\partial}

\newcommand{\cil}{\xrightarrow{\mcl}} 
\newcommand{\asc}{\xrightarrow{a.s.}} 
\newcommand{\argmax}{\mathop{\rm argmax}}
\newcommand{\diag}{\mathop{\rm diag}} 
\newcommand{\mbbi}{\mathbbm{1}} 

\def\ds#1{\displaystyle{#1}}
\def\nn{\nonumber}

\def\sumj{\sum_{j=1}^{n}}

\def\dim{\mathrm{dim}}

\def\tes{\hat{\theta}_{n}}
\def\aes{\hat{\alpha}_{n}}

\def\ges{\hat{\gamma}_{n}}

\title[]
{Schwartz type model selection for ergodic stochastic differential equation models}
\date{\today}

\author{Shoichi Eguchi and Yuma Uehara}

\address[Yuma Uehara, corresponding author]{The Institute of Statistical Mathematics, Japan, 10-3 Midori-cho, Tachikawa, Tokyo 190-8562, Japan}
\email{y-uehara@ism.ac.jp}

\begin{document}

\maketitle

\begin{abstract}
We study theoretical foundation of model comparison for ergodic stochastic differential equation (SDE) models and an extension of the applicable scope of the conventional Bayesian information criterion.
Different from previous studies, we suppose that the candidate models are possibly misspecified models, and we consider both Wiener and a pure-jump L\'{e}vy noise driven SDE.
Based on the asymptotic behavior of the marginal quasi-log likelihood, the Schwarz type statistics and stepwise model selection procedure are proposed.
We also prove the model selection consistency of the proposed statistics with respect to an optimal model.
We conduct some numerical experiments and they support our theoretical findings.
\end{abstract}

\section{Introduction} 

We suppose that the data-generating process $X$ is defined on the stochastic basis $(\Omega, \mcf, \mcf_t, P)$ and it is the solution of the one-dimensional stochastic differential equation written as:
\begin{equation}\label{yu:tmodel}
dX_t=A(X_t)dt+C(X_{t-})dZ_t,
\end{equation}
where:
\begin{itemize}
\item The coefficients $A$ and $C$ are Lipschitz continuous.
\item The driving noise $Z$ is a standard Wiener process or pure-jump L\'{e}vy process satisfying that for any $q>0$,
\begin{equation}\label{yu:momcon}
E[Z_1]=0, \quad E[Z_1^2]=1, \quad E[|Z_1|^q]<\infty.
\end{equation}
\item The initial variable $X_0$ is independent of $Z$, and 
\begin{equation*}
\mcf_t=\sig(X_0)\vee\sig(Z_s|s\leq t)
\end{equation*}
\end{itemize}
As the observations from $X$, we consider the discrete but high-frequency samples $(X_{t_j^n})_{j=0}^n$ with 
\begin{equation*}
t_j^n:=jh_n, \quad T_n:=nh_n\to\infty, \quad nh_n^2\to0.
\end{equation*}
For $(X_{t_j^n})_{j=0}^n$, $M_1\times M_2$ candidate models are supposed to be given.
Here, for each $m_{1}\in\{1,\dots, M_{1}\}$ and $m_{2}\in\{1,\dots, M_{2}\}$, the candidate model $\mcm_{m_1,m_2}$ is expressed as:
\begin{align}\label{yu:canmodel}
dX_t=a_{m_{2}}(X_t,\al_{m_{2}})dt+c_{m_{1}}(X_{t-},\gam_{m_{1}})dZ_t,
\end{align}
and the functional form of $(c_{m_1}(\cdot,\cdot),a_{m_2}(\cdot,\cdot))$ is known except for the $p_{\gam_{m_1}}$ and $p_{\al_{m_2}}$-dimensional unknown parameters $\gamma_{m_{1}}$ and $\alpha_{m_{2}}$ being elements of the bounded convex domains $\Theta_{\gamma_{m_{1}}}\subset\mathbb{R}^{p_{\gamma_{m_{1}}}}$ and $\Theta_{\alpha_{m_{2}}}\subset\mathbb{R}^{p_{\alpha_{m_{2}}}}$.
The main objective of this paper is to give a model selection procedure for extracting an ``optimal" model $\mcm_{m_1^\star,m_2^\star}$ among $\mcm:=\{\mcm_{m_1,m_2}| m_{1}\in\{1,\dots, M_{1}\}, m_{2}\in\{1,\dots, M_{2}\}\}$ which reflects the feature of $X$ well.

For selecting an appropriate model from the data in hand quantitively, information criteria are one of the most convenient and powerful tools, and have widely been used in many fields.
Their origin dates back to Akaike information criterion (AIC) introduced in \cite{Aka73,Aka74} which puts an importance on prediction, and after that, various kinds of criteria have been produced up to the present, for their comprehensive overview, see \cite{BurAnd02}, \cite{ClaHjo08}, and \cite{KonKit08}.
Among them, this paper especially sheds light on Bayesian information criterion (BIC) introduced by \cite{Sch78}.
It is based on an approximation up to $O_p(1)$-term of log-marginal likelihood, and its original form is as follows:
\begin{equation}\label{BIC}
\text{BIC}_n=-2l_n(\tes^{\text{MLE}})+p\log n,
\end{equation}
where $l_n$, $\tes^{\text{MLE}}$, and $p$ stand for the log-likelihood function, maximum likelihood estimator, and dimension of the parameter including the subject model.
However, since the closed form of the transition density of $X$ is unavailable in general, to conduct some feasible statistical analysis, we cannot rely on its genuine likelihood; this implies that the conventional likelihood based (Bayesian) information criteria are unpractical in our setting.
Such a problem often occurs when discrete time observations are obtained from a continuous time process, and to avoid it, the replacement of a genuine likelihood by some quasi-likelihood is effective not only for estimating parameters included in a subject model but also for constructing (quasi-)information criteria, for instance, see \cite{Uch10}, \cite{FujUch14}, \cite{EguMas18a} (ergodic diffusion model), \cite{UchYos16} (stochastic regression model), and \cite{FasKim17} (CARMA process). 
Especially, \cite{EguMas18a} used the Gaussian quasi-likelihood in place of the genuine likelihood, and derived quasi-Bayesian information criterion (QBIC) under the conditions: the driving noise is a standard Wiener process, and for each candidate model, there exist $\gam_{m_1,0}\in\Theta_{m_1}$ and $\al_{m_2,0}\in\Theta_{m_2}$ satisfying $c_{m_1}(x,\gam_{m_1,0})\equiv C(x)$ and $a_{m_2}(x,\al_{m_2,0})\equiv A(x)$, respectively.
Moreover, by using the difference of the small time activity of the drift and diffusion terms, the paper also gave the two-step QBIC which selects each term separately, and reduces the computational load.
In the paper, the model selection consistency of the QBIC is shown only for nested case.
In such a case, by considering the (largest) model which contains all candidate models, regularized methods can also be used in the same purpose. Concerning the regularized method for SDE models, for example, see \cite{GreIac12} and \cite{MasShi17}.

While when it comes to the estimation of the parameters $\gam_{m_1}$ and $\al_{m_2}$, the Gaussian quasi maximum likelihood estimator (GQMLE) works well for a much broader situation: the driving noise is a standard Wiener process or pure-jump L\'{e}vy process with \eqref{yu:momcon}, and either or both of the drift and scale coefficients are misspecified.
For the technical account of the GQMLE for ergodic SDE models, see \cite{Yos92}, \cite{Kes97}, \cite{UchYos11}, \cite{UchYos12}, \cite{Mas13-1}, and \cite{Ueh18}.
These results naturally provides us an insight that the aforementioned QBIC is also theoretically valid for the broader situation, and has the model selection consistency even if a non-nested model is contained in candidate models.
In this paper, we will show that the insight is true.
More specifically, we will give the QBIC building on the stochastic expansion of the log-marginal Gaussian quasi-likelihood.
Although the convergence rate of the GQMLE differs in the L\'{e}vy driven or misspecified case, the form is the same as the correctly specified diffusion case, that is, a unified model selection criteria for ergodic SDE models is provided.
We will also show the model selection consistency of the QBIC.

The rest of this paper is as follows: Section \ref{sec_ass} provides the notations and assumptions used throughout this paper. In Section \ref{sec_res}, the main result of this paper is given. Section \ref{sec_sim} exhibits some numerical experiments. The technical proofs of the main results are summarized in Appendix.



\section{Notations and Assumptions}\label{sec_ass}
For notational convenience, we previously introduce some symbols used in the rest of this paper.
\begin{itemize} 
\item For any vector $x$, $x^{(j)}$ represents $j$-th element of $x$, and we write $x^{\otimes2}=x^\top x$ where $^\top$ denotes the transpose operator.
\item $\p_x$ is referred to as a differential operator with respect to any variable $x$.
\item $x_n\lesssim y_n$ implies that there exists a positive constant $C$ being independent of $n$ satisfying $x_n\leq Cy_n$ for all large enough $n$.
\item For any set $S$, $\bar{S}$ denotes its closure.
\item We write $Y_j=Y_{t_j}$ and $\D_j Y:=Y_j-Y_{j-1}$ for any stochastic process $(Y_t)_{t\in\mbbr^+}$.
\item For any matrix valued function $f$ on $\mbbr\times\Theta$, we write $f_s(\theta)=f(X_s,\theta)$; especially we write $f_j(\theta)=f(X_j,\theta)$.
\item $I_p$ represents the $p$-dimensional identity matrix. 
\item $\nu_0$ represents the L\'{e}vy measure of $Z$.
\item $P_t(x,\cdot)$ denotes the transition probability of $X$. 
\item Given a function $\rho:\mathbb{R}\to\mathbb{R}^+$ and a signed measure $m$ on a one-dimensional Borel space, we write
\begin{equation}\nn
||m||_\rho=\sup\left\{|m(f)|:\mbox{$f$ is $\mathbb{R}$-valued, $m$-measurable and satisfies $|f|\leq\rho$}\right\}.
\end{equation}
\item $\mca$ and $\tilde{\mca}$ stand for the infinitesimal generator and extended generator of $X$, respectively.
\end{itemize}

In the next section, we will first give the stochastic expansion of the log-marginal Gaussian quasi-likelihood for the following model:
\begin{equation}\label{ten:model}
dX_t=a(X_t,\al)dt+c(X_{t-},\gam)dZ_t,
\end{equation}
where similar to \eqref{yu:canmodel}, the coefficients have the unknown $p_\gam$-dimensional parameter $\gam$ and $p_\al$-dimensional parameter $\al$.
They are supposed to be elements of bounded convex domains $\Theta_\gam$ and $\Theta_\al$, and for the sake of convenience, we write $\theta=(\gam,\al)$ and $\Theta_\gam\times \Theta_\al:=\Theta$.
We also assume that either or both of the drift and scale coefficients are possibly misspecified.
Especially, we say that the model setting is semi-misspecified diffusion case when the driving noise is a standard Wiener process, the scale coefficient is correctly specified and the drift coefficient is misspecified.

Below, we table the assumptions for our main result.

\begin{Assumption}\label{Moments}
$Z$ is a standard Wiener process, or a pure-jump L\'{e}vy process satisfying that: 
$E[Z_1]=0$, $E[Z_1^2]=1$, and $E[|Z_1|^q]<\infty$ for all $q>0$.
Furthermore, the Blumenthal-Getoor index (BG-index) of $Z$ is smaller than 2, that is, 
\begin{equation*}
\beta:=\inf_\gam\left\{\gam\geq0: \int_{|z|\leq1}|z|^\gam\nu_0(dz)<\infty\right\}<2.
\end{equation*}
\end{Assumption}

\begin{Assumption}\label{Stability}
\begin{enumerate}
\item
There exists a probability measure $\pi_0$ such that for every $q>0$, we can find constants $a>0$ and $C_q>0$ for which 
\begin{equation}\label{Ergodicity}
\sup_{t\in\mathbb{R}_{+}} \exp(at) ||P_t(x,\cdot)-\pi_0(\cdot)||_{h_q} \leq C_qh_q(x),
\end{equation}
for any $x\in\mathbb{R}$ where $h_q(x):=1+|x|^q$.
\item 
For any $q>0$,  we have
\begin{equation}
\sup_{t\in\mathbb{R}_{+}}E[|X_t|^q]<\infty. 
\end{equation}
\end{enumerate}
\end{Assumption}

Let $\pi_1$ and $\pi_2$ be the prior densities for $\gam$ and $\al$, respectively.
\begin{Assumption}\label{Prior}
The prior densities $\pi_1$ and $\pi_2$ are continuous, and fullfil that
\begin{equation*}
\sup_{\gam\in\Theta_\gam} \pi_1(\gam) \vee \sup_{\al\in\Theta_\al} \pi_2(\al) <\infty.
\end{equation*}
\end{Assumption}

We define an {\it optimal} value $\theta^\star:=(\gam^\star,\al^\star)$ of $\theta$ being chosen arbitrary from the sets $\displaystyle\argmax_{\gam\in\bar{\Theta}_\gam}\mbbg_1(\gam)$ and $\displaystyle\argmax_{\al\in\bar{\Theta}_\al}\mbbg_2(\al)$
for $\mbbr$-valued functions $\mbbg_1(\cdot)$ (resp. $\mbbg_{2}(\cdot)$) on $\Theta_\gam$ (resp. $\Theta_\al$) defined by
\begin{align}
&\mbbg_1(\gam):=-\int_\mbbr \left(\log c^2(x,\gam)+\frac{C^2(x)}{c^2(x,\gam)}\right)\pi_0(dx), \label{rf:con.scale}\\
&\mbbg_2(\al):=-\int_\mbbr c(x,\gam^\star)^{-2}(A(x)-a(x,\al))^2\pi_0(dx). \label{rf:con.drift}
\end{align}
From the expression of $\mbbg_1$, $\gam^\star$ can be regarded as an element in $\Theta_\gam$ minimizing the Stein's loss; $\al^\star$ as an element in $\Theta_\al$ minimizing $L_2$-loss.
Recall that $\Theta=\Theta_\gam\times\Theta_\al$ is supposed to be a bounded convex domain.
Then, we assume that:
\begin{Assumption}\label{Identifiability}
\begin{itemize}
\item $\theta^\star$ is unique and is in $\Theta$.
\item There exist positive constants $\chi_\gam$ and $\chi_\al$ such that for all $(\gam,\al)\in\Theta$,
\begin{align}
&\mbbg_1(\gam)-\mbbg_1(\gam^\star)\leq-\chi_\gam|\gam-\gam^\star|^2,\\
&\mbbg_2(\al)-\mbbg_2(\al^\star)\leq-\chi_\al|\al-\al^\star|^2.
\end{align}
\end{itemize}
\end{Assumption}

\begin{Assumption}\label{Smoothness}
\begin{enumerate}
\item The coefficients $A$ and $C$ are Lipschitz continuous and twice differentiable, and their first and second derivatives are of at most polynomial growth.
\item The drift coefficient $a(\cdot,\al^\star)$ and scale coefficient $c(\cdot,\gam^\star)$ are Lipschitz continuous, and $c(x,\gam)\neq0$ for every $(x,\gam)$.
\item For each $i \in \left\{0,1\right\}$ and $k \in \left\{0,\dots,5\right\}$, the following conditions hold:
\begin{itemize}
\item The coefficients $a$ and $c$ admit extension in $\mathcal{C}(\mathbb{R}\times\bar{\Theta})$ and have the partial derivatives $(\partial_x^i \partial_\alpha^k a, \partial_x^i \partial_\gamma^k c)$ possessing extension in $\mathcal{C}(\mathbb{R}\times\bar{\Theta})$.
\item There exists nonnegative constant $C_{(i,k)}$ satisfying
\begin{equation}\label{polynomial}
\sup_{(x,\alpha,\gamma) \in \mathbb{R} \times \bar{\Theta}_\alpha \times \bar{\Theta}_\gamma}\frac{1}{1+|x|^{C_{(i,k)}}}\left\{|\partial_x^i\partial_\alpha^ka(x,\alpha)|+|\partial_x^i\partial_\gamma^kc(x,\gamma)|+|c^{-1}(x,\gamma)|\right\}<\infty.
\end{equation}
\end{itemize}
\end{enumerate}
\end{Assumption}

Define the $p_\gam\times p_\gam$-matrix $\mci_\gam$ and $p_\al\times p_\al$-matrix $\mci_\al$ by:
\begin{align}
&\mci_\gam=4\int_\mbbr \frac{(\p_\gam c(x,\gam^\star))^{\otimes2}}{c^4(x,\gam^\star)}C^2(x)\pi_0(dx)-2\int_\mbbr\frac{\p_\gam^{\otimes2}c(x,\gam^\star)c(x,\gam^\star)-(\p_\gam c(x,\gam^\star))^{\otimes2}}{c^4(x,\gam^\star)}(C^2(x)-c^2(x,\gam^\star))\pi_0(dx),\label{yu:fish1}\\
&\mci_\al=2\int_\mbbr\frac{(\p_\al a(x,\al^\star))^{\otimes2}}{c^2(x,\gam^\star)}\pi_0(dx)-2\int_\mbbr\frac{\p_\al^{\otimes2} a(x,\al^\star)}{c^2(x,\gam^\star)}(A(x)-a(x,\al^\star))\pi_0(dx)\label{yu:fish2}.
\end{align}
In the correctly specified case, since under Assumption \ref{Identifiability}, we have
\begin{equation*}
c(x,\gam^\star)= C(x),\quad a(x,\al^\star)=A(x), \qquad \pi_0-a.s.,
\end{equation*}
$\mci_\gam$ and $\mci_\al$ are reduced to
\begin{align*}
&\mci_\gam=4\int_\mbbr \frac{(\p_\gam c(x,\gam^\star))^{\otimes2}}{c^2(x,\gam^\star)}\pi_0(dx),\\
&\mci_\al=2\int_\mbbr\frac{(\p_\al a(x,\al^\star))^{\otimes2}}{c^2(x,\gam^\star)}\pi_0(dx).
\end{align*}

\begin{Assumption}\label{Fisher}
$\mci_\gam$ and $\mci_\al$ are positive definite.
\end{Assumption}

In the rest of this section, we give a brief overview of the stepwise Gaussian quasi-likelihood method for \eqref{ten:model}, and introduce some theoretical results under Assumptions \ref{Moments}-\ref{Fisher}.
We consider the following stepwise Gaussian quasi-likelihood (GQL) functions $\mbbg_{1,n}$ and $\mbbg_{2,n}$ on $\Theta_\gam$ and $\Theta_\al$: 
\begin{align}
&\mbbg_{1,n}(\gam)=-\frac{1}{h_n}\sumj \left\{h_n\log c^2_{j-1}(\gam)+\frac{(\D_j X)^2}{c^2_{j-1}(\gam)}\right\} \label{rf:scale},\\
&\mbbg_{2,n}(\al)=-\sumj \frac{(\D_j X-h_na_{j-1}(\al))^2}{h_nc^2_{j-1}(\ges)} \label{rf:drift}
\end{align}
For such functions, we define the (stepwise) Gaussian quasi maximum likelihood estimator (GQMLE) $\tes:=(\ges,\aes)$ in the following manner:
\begin{align*}
&\ges\in\argmax_{\gam\in\bar{\Theta}_\gam}\mbbg_{1,n}(\gam),
\nonumber\\
&\aes\in\argmax_{\al\in\bar{\Theta}_\al}\mbbg_{2,n}(\al).
\nonumber
\end{align*}

\begin{Rem}\label{jointstepwise}
Different from the low frequently observed case, the stepwise type estimator exhibits the same asymptotics as the joint type estimator defined as:
\begin{equation*}
\tilde{\theta}_n\in\argmax_{\theta\in\bar{\Theta}} \left[-\frac{1}{h_n}\sumj \left\{h_n\log c^2_{j-1}(\gam)+\frac{(\D_j X-h_na_{j-1}(\al))^2}{c^2_{j-1}(\gam)}\right\}\right],
\end{equation*}
while possessing the stability of optimization for calculating the estimator.
This is because the small time behavior is dominated by the scale term, and thus theoretically, $h_na_{j-1}(\al)$ does not affect the estimation of $\gam$.
\end{Rem}

By making use of the estimates in the papers \cite{Yos92}, \cite{Kes97}, \cite{UchYos11}, \cite{UchYos12}, \cite{Mas13-1}, and \cite{Ueh18}, the following asymptotic results about $\tes$ can be obtained (or directly follows) under Assumption \ref{Moments}-\ref{Fisher}: let $A_n:=\diag\{a_n I_{p_\gam}, \sqrt{T_n} I_{p_\al}\}$ where $a_n=\sqrt{n}$ in the correctly specified or semi-misspecified diffusion case, and otherwise, $a_n=\sqrt{T_n}$. 
Then we have
\begin{itemize}
\item Tail probability estimates:
for any $L>0$ and $r>0$, there exists a positive constant $C_L$ such that
 \begin{equation}\label{eq: TPE}
\sup_{n\in\mbbn} P\left(\left|A_n(\tes-\theta^\star)\right|>r\right)\leq \frac{C_L}{r^L}.
 \end{equation}
\item Asymptotic normality:
\begin{equation*}
A_n(\tes-\theta^{\star})\cil N(0, \mci^{-1}\Sig (\mci^{-1})^\top),
\end{equation*}
where $\mci=\begin{pmatrix}\mci_\gam & O \\ \mci_{\al\gam} & \mci_{\al}\end{pmatrix}$ with
\begin{equation*}
\mci_{\al\gam}=2\int_\mbbr \p_\al a(x,\al^\star)\p_\gam^\top c^{-2}(x,\gam^\star)(a(x,\al^\star)-A(x))  \pi_0(dx),
\end{equation*} 
and the form of $\Sig:=\begin{pmatrix}\Sig_\gam&\Sig_{\al\gam}\\\Sig_{\al\gam}^\top&\Sig_{\al}\end{pmatrix}$ is given as follows: 
\begin{enumerate}
\item Correctly specified diffusion case:
\begin{equation*}
\Sig=2\mci=2\diag\{\mci_\gam, \mci_\al\}=\begin{pmatrix}8\int_\mbbr \frac{(\p_\gam c(x,\gam^\star))^{\otimes2}}{c^2(x,\gam^\star)}\pi_0(dx) & O\\ O &4\int_\mbbr\frac{(\p_\al a(x,\al^\star))^{\otimes2}}{c^2(x,\gam^\star)}\pi_0(dx) \end{pmatrix},
\end{equation*}
hence the asymptotic variance can be simply written as
\begin{equation*}
\mci^{-1}\Sig (\mci^{-1})^\top=2 \mci^{-1}=\begin{pmatrix}\frac{1}{2}\int_\mbbr \frac{(\p_\gam c(x,\gam^\star))^{\otimes2}}{c^2(x,\gam^\star)}\pi_0(dx) & O\\ O &\int_\mbbr\frac{(\p_\al a(x,\al^\star))^{\otimes2}}{c^2(x,\gam^\star)}\pi_0(dx) \end{pmatrix}.
\end{equation*}
\item Semi-misspecified diffusion case:
\begin{align*}
&\Sig_\gam =8\int_\mbbr \frac{(\p_\gam c(x,\gam^\star))^{\otimes2}}{c^2(x,\gam^\star)}\pi_0(dx) \\ 
&\Sig_{\al\gam}=0\\
&\Sig_\al=4\int \left[\left(\frac{\p_\al a(x,\al^\star)}{c^2(x,\gam^\star)}-\p_x f(x)\right)C(x)\right]^{\otimes 2}\pi_0(dx)
\end{align*}
where the function $f$ is the solution of the following Poisson equation:
\begin{align*}
\mca f^{(j)}(x)&=\frac{\p_{\al^{(j)}} a(x,\al^\star)}{c^2(x,\gam^\star)}(A(x)-a(x,\al^\star)),
\end{align*}
for $j\in\{1,\dots, p_\al\}$.

\item Misspecified diffusion case:
\begin{align*}
&\Sig_\gam =4\int (\p_x f_1(x) C(x))^{\otimes 2}\pi_0(dx)\\ 
&\Sig_{\al\gam}=4\int \left(\frac{\p_\al a(x,\al^\star)}{c^2(x,\gam^\star)}-\p_x f_2(x)\right)C^2(x)(\p_x f_1(x))^\top \pi_0(dx)\\
&\Sig_\al=4\int \left[\left(\frac{\p_\al a(x,\al^\star)}{c^2(x,\gam^\star)}-\p_x f_2(x)\right)C(x)\right]^{\otimes 2}\pi_0(dx)
\end{align*}
where the functions $f_1$ and $f_2$ are the solution of the following Poisson equations:
\begin{align*}
\mca f_1^{(j_1)}(x)&=\frac{\p_{\gam^{(j_1)}} c(x,\gam^\star)}{c^3(x,\gam^\star)}(c^2(x,\gam^\star)-C^2(x)), \\
\mca f_2^{(j_2)}(x)&=\frac{\p_{\al^{(j_2)}} a(x,\al^\star)}{c^2(x,\gam^\star)}(A(x)-a(x,\al^\star)),
\end{align*}
for $j_1\in\{1,\dots, p_\gam\}$ and $j_2\in\{1,\dots, p_\al\}$.
\item L\'{e}vy driven case (both correctly specified and misspecified case):
\begin{align*}
&\Sig_\gam=4\int_\mbbr\int_\mbbr\left(\frac{\p_\gam c(x,\gam^\star)}{c^3(x,\gam^\star)}C^2(x)z^2+g_1(x+C(x)z)-g_1(x)\right)^{\otimes2}\pi_0(dx)\nu_0(dz),\\
&\Sig_{\al\gam}=-4\int_\mbbr\int_\mbbr\left(\frac{\p_\gam c(x,\gam^\star)}{c^3(x,\gam^\star)}C^2(x)z^2+g_1(x+C(x)z)-g_1(x)\right)\\
&\qquad \qquad \quad\left(\frac{\p_\al a(x,\al^\star)}{c^2(x,\gam^\star)}C(x)z+g_2(x+C(x)z)-g_2(x)\right)^\top\pi_0(dx)\nu_0(dz),\\
&\Sig_\al=4\int_\mbbr\int_\mbbr\left(\frac{\p_\al a(x,\al^\star)}{c^2(x,\gam^\star)}C(x)z+g_2(x+C(x)z)-g_2(x)\right)^{\otimes2}\pi_0(dx)\nu_0(dz),
\end{align*}
where the functions $g_1$ and $g_2$ are the solution of the following extended Poisson equations:
\begin{align*}
\tilde{\mca} g_1^{(j_1)}(x)&=-\frac{\p_{\gam^{(j_1)}} c(x,\gam^\star)}{c^3(x,\gam^\star)}(c^2(x,\gam^\star)-C^2(x)), \\
\tilde{\mca} g_2^{(j_2)}(x)&=-\frac{\p_{\al^{(j_2)}} a(x,\al^\star)}{c^2(x,\gam^\star)}(A(x)-a(x,\al^\star)),
\end{align*}
for $j_1\in\{1,\dots, p_\gam\}$ and $j_2\in\{1,\dots, p_\al\}$ (In the correctly specified case, $g_1$ and $g_2$ are identically 0).
\end{enumerate}
\end{itemize}

We note that in the (semi-)misspecified diffusion case, the asymptotic results on the stepwise GQMLE is not verified. However, by taking the same route to \cite[Theorem 3.1]{Ueh18}, we can easily show the tail probability estimates. Concerning asymptotic normality, it can also be derived from the argument of Remark \ref{jointstepwise}.

\begin{Rem}
The theory of Poisson equation and extended Poisson equation plays an important role for dealing with the misspecification effect.
The former equation corresponds the generator of diffusion processes, and the latter one does the extended generator of Feller Markov processes.
The existence and regularity conditions for their solutions are discussed in \cite{ParVer01}, \cite{VerKul11}, and \cite{Ueh18}, and in the diffusion case, \cite[Remark 2.2]{UchYos11} provides the explicit form of $g_1$ and $g_2$ under $p_\gam=p_\al=1$.
\end{Rem}

\section{Main results}\label{sec_res}

Building on the stepwise Gaussian quasi-likelihood $\mbbg_{1,n}$ and $\mbbg_{2,n}$, the next theorem gives the stochastic expansion of the log-marginal quasi-likelihood which is the main result of this paper:
\begin{Thm}\label{YU:se}
If Assumptions \ref{Moments}-\ref{Fisher} are satisfied for the statistical model \eqref{ten:model}, we have
\begin{align*}
&\log\left(\int_{\Theta_\gam}\exp\left(\mbbg_{1,n}(\gam)\right)\pi_1(\gam)d\gam\right)=\mbbg_{1,n}(\ges)-\frac{1}{2}p_\gam \log n+\log\pi_1\left(\gam^\star\right)+\frac{p_\gam}{2}\log 2\pi-\frac{1}{2}\log \det \mci_\gam+o_p\left(1\right),\\
&\log\left(\int_{\Theta_\al}\exp\left(\mbbg_{2,n}(\al)\right)\pi_2(\al)d\al\right)=\mbbg_{2,n}(\aes)-\frac{1}{2}p_\al \log T_n+\log \pi_2(\al^\star)+\frac{p_\al}{2}\log 2\pi-\frac{1}{2}\log\det \mci_\al+o_p(1).
\end{align*}
\end{Thm}

In the present settings, although the scale estimator has the two kinds of convergence rates depending on the model setups, Theorem \ref{YU:se} holds regardless of the convergence rate of the scale estimator.
By ignoring the $O_p(1)$ terms in each expansion, we define the two-step quasi-Bayesian information criteria (QBIC) by
\begin{align*}
&\mbox{QBIC}_{1,n}=\mbbg_{1,n}(\ges)-\frac{1}{2}p_\gam \log n,\\
&\mbox{QBIC}_{2,n}=\mbbg_{2,n}(\aes)-\frac{1}{2}p_\al \log(T_n).
\end{align*}
 

Next, we consider model selection consistency of the proposed information criteria.
Suppose that candidates for drift and scale coefficients are given as
\begin{align}
& c_{1}(x,\gamma_{1}),\ldots,c_{M_{1}}(x,\gamma_{M_{1}}), \label{se:ms.c} \\
& a_{1}(x,\alpha_{1}),\dots,a_{M_{2}}(x,\alpha_{M_{2}}), \label{se:ms.a} 
\end{align}
where $\gamma_{m_{1}}\in\Theta_{\gamma_{m_{1}}}\subset\mathbb{R}^{p_{\gamma_{m_{1}}}}$ for any $m_{1}\leq M_{1}$ and $\alpha_{m_{2}}\in\Theta_{\alpha_{m_{2}}}\subset\mathbb{R}^{p_{\alpha_{m_{2}}}}$ for any $m_{2}\leq M_{2}$.
Then, each candidate model $\mathcal{M}_{m_{1},m_{2}}$ is given by
\begin{align*}
dX_t=a_{m_{2}}(X_t,\al_{m_{2}})dt+c_{m_{1}}(X_{t-},\gam_{m_{1}})dZ_t.
\end{align*}
In each candidate model $\mathcal{M}_{m_{1},m_{2}}$, the functions \eqref{rf:scale} and \eqref{rf:con.scale} are denoted by $G_{1,n}^{(m_{1})}$ and $G_{1}^{(m_{1})}$, respectively.
The functions $G_{2,n}^{(m_{2}|m_{1})}$ and $G_{2}^{(m_{2}|m_{1})}$ correspond to \eqref{rf:drift} and \eqref{rf:con.drift} with $\gamma_{m_{1}}$. 
Using the QBIC, we propose the stepwise model selection as follows. 
\begin{itemize}
\item[(i)] We select the best scale coefficient $c_{\hat{m}_{1,n}}$ among \eqref{se:ms.c}, where $\hat{m}_{1,n}$ satisfies $\{\hat{m}_{1,n}\}=\argmax_{m_{1}}\mathrm{QBIC}_{1,n}^{(m_{1})}$ with
\begin{align*}
&\mathrm{QBIC}_{1,n}^{(m_{1})}=G_{1,n}^{(m_{1})}(\hat{\gamma}_{m_{1},n})-p_{\gamma_{m_{1}}}\log n,\\
&\hat{\gamma}_{m_{1},n}\in\argmax_{\gam_{m_{1}}\in\bar{\Theta}_{\gam_{m_{1}}}}\mbbg_{1,n}^{(m_{1})}(\gam_{m_{1}}).
\end{align*}

\item[(ii)] Under $c_{\hat{m}_{1,n}}$ and $\hat{\gamma}_{\hat{m}_{1,n},n}$, we select the best drift coefficient with index $\hat{m}_{2,n}$ such that $\{\hat{m}_{2,n}\}=\argmax_{m_{2}}\mathrm{QBIC}_{2,n}^{(m_{2}|\hat{m}_{1,n})}$, where
\begin{align*}
&\mathrm{QBIC}_{2,n}^{(m_{2})}=G_{2,n}^{(m_{2}|\hat{m}_{1,n})}(\hat{\alpha}_{m_{2},n})-p_{\alpha_{m_{2}}}\log(T_{n}),\\
&\hat{\alpha}_{m_{2},n}\in\argmax_{\alpha_{m_{2}}\in\bar{\Theta}_{\alpha_{m_{2}}}}\mbbg_{2,n}^{(m_{2}|\hat{m}_{1,n})}(\alpha_{m_{2}}).
\end{align*}

\end{itemize} 
Through this procedure, we can obtain the model $\mcm_{\hat{m}_{1,n},\hat{m}_{2,n}}$ as the final best model among the candidates described by \eqref{se:ms.c} and \eqref{se:ms.a}.

The {\it optimal value} $\theta_{m_{1},m_{2}}^{\star}=(\alpha_{m_{2}}^{\star},\gamma_{m_{1}}^{\star})$ of $\mathcal{M}_{m_{1},m_{2}}$ is defined in a similar manner as the previous section.
We assume that the model indexes $m_{1}^{\star}$ and $m_{2}^{\star}$ are uniquely given as follows: 
\begin{align*}
\{m_{1}^{\star}\}&= \argmin_{m_{1}\in\mathfrak{M}_{1}}\dim(\Theta_{\gamma_{m_{1}}}),\\
\{m_{2}^{\star}\}&= \argmin_{m_{2}\in\mathfrak{M}_{2}}\dim(\Theta_{\alpha_{m_{2}}}),
\end{align*}
where $\mathfrak{M}_{1}=\argmax_{1\leq m_{1}\leq M_{1}}G_{1}^{(m_{1})}(\gamma_{m_{1}}^{\star})$ and $\mathfrak{M}_{2}=\argmax_{1\leq m_{2}\leq M_{2}}G_{2}^{(m_{2}|m_{1}^{\star})}(\alpha_{m_{2}}^{\star})$.
Then, we say that $\mathcal{M}_{m_{1}^{\star},m_{2}^{\star}}$ is the {\it optimal model}.
That is, the optimal model consists of the elements of optimal model sets $\mathfrak{M}_{1}$ and $\mathfrak{M}_{2}$ which have the smallest dimension.
The following theorem means that the proposed criteria and model selection method have the model selection consistency.

\begin{Thm} \label{thm:mod.cons}
Suppose that Assumptions \ref{Moments}-\ref{Fisher} hold for the all candidate models and that these exists a $\mathcal{M}_{m_{1}^{\star},m_{2}^{\star}}$ is the optimal model.
Let $m_{1}\in\{1,\ldots,M_{1}\}\backslash\{m_{1}^{\star}\}$ and $m_{2}\in\{1,\ldots,M_{2}\}\backslash\{m_{2}^{\star}\}$. 
Then the model selection consistency of the proposed QBIC hold in the following senses.
\begin{align*}
& \lim_{n\to\infty}\mathbb{P}\left(\mathrm{QBIC}_{1,n}^{(m_{1}^{\star})}-\mathrm{QBIC}_{1,n}^{(m_{1})}>0\right) =1, \\
& \lim_{n\to\infty}\mathbb{P}\left(\mathrm{QBIC}_{2,n}^{(m_{2}^{\star}|\hat{m}_{1,n})}-\mathrm{QBIC}_{2,n}^{(m_{2}|\hat{m}_{1,n})}>0\right)=1.
\end{align*} 
\label{se:thm.modcon}
\end{Thm}

\begin{Rem}
We here consider the case where there are several optimal models.
Then, we define the optimal model index sets $\mathfrak{M}_{1}^{\star}$ and $\mathfrak{M}_{2}^{\star}$ by
\begin{align*}
\mathfrak{M}_{1}^{\star}&= \argmin_{m_{1}\in\mathfrak{M}_{1}}\dim(\Theta_{\gamma_{m_{1}}}),\\
\mathfrak{M}_{2}^{\star}&= \argmin_{m_{2}\in\mathfrak{M}_{2}}\dim(\Theta_{\alpha_{m_{2}}}),
\end{align*}
respectively.
Applying the proof of Theorem \ref{se:thm.modcon} for each elements of $\mathfrak{M}_{1}^{\star}$ and $\mathfrak{M}_{2}^{\star}$, we can show the model selection consistency with respect to the optimal model sets.
\end{Rem}

\section{Numerical experiments}\label{sec_sim}

In this section, we present simulation results to observe finite-sample performance of the proposed QBIC.
We use {\tt yuima} package on R (see \cite{YUIMA14}) for generating data.
In the examples below,  all the Monte Carlo trials are based on 1000 independent sample paths, and the simulations are done for $(h_{n},T_{n})=(0.01,10), (0.005,10), (0.01,50)$, and $(0.005,50)$ (hence in each case, $n=1000, 2000, 5000$, and $10000$).
We simulate the model selection frequencies by using proposed QBIC and compute the model weight $w_{m_{1},m_{2}}$ (\cite[Section 6.4.5]{BurAnd02}) defined by
\begin{align}
\begin{split}
w_{m_{1},m_{2}}&=\frac{\ds{\exp\Big\{-\frac{1}{2}\big(\mathrm{QBIC}_{1,n}^{(m_{1})}-\mathrm{QBIC}_{1,n}^{(\hat{m}_{1,n})}\big)\Big\}}}{\ds{\sum_{k=1}^{M_{1}}\exp\Big\{-\frac{1}{2}\big(\mathrm{QBIC}_{1,n}^{(k)}-\mathrm{QBIC}_{1,n}^{(\hat{m}_{1,n})}\big)\Big\}}} \\
&\quad\times\frac{\ds{\exp\Big\{-\frac{1}{2}\big(\mathrm{QBIC}_{2,n}^{(m_{2}|m_{1})}-\mathrm{QBIC}_{n}^{(\hat{m}_{2,n}|m_{1})}\big)\Big\}}}{\ds{\sum_{\ell=1}^{M_{2}}\exp\Big\{-\frac{1}{2}\big(\mathrm{QBIC}_{2,n}^{(\ell|m_{1})}-\mathrm{QBIC}_{2,n}^{(\hat{m}_{2,n}|m_{1})}\big)\Big\}}}\times100. 
\label{def:weight} 
\end{split}
\end{align}
The model weight can be used to empirically quantify relative frequency (percentage) of the model selection from a single data set.
The model which has the highest $w_{m_{1},m_{2}}$ value is regarded as the most probable model. 
From its definition \eqref{def:weight}, $w_{m_{1},m_{2}}$ satisfies the equation $\sum_{k=1}^{M_{1}}\sum_{\ell=1}^{M_{2}}w_{k,\ell}=100$.

\subsection{Ergodic diffusion model}\label{subsec_sim1}

Suppose that we have a sample $\mbX_{n}=(X_{t_{j}})_{j=0}^{n}$ with $t_{j}=jh_{n}$ from the true model
\begin{align*}
dX_{t}=-\frac{1}{2}X_{t}dt+dw_{t},\quad t\in[0,T_{n}],\quad X_{0}=0,
\end{align*}
where $T_{n}=nh_{n}$, and $w$ is a one-dimensional standard Wiener process. 
We consider the following scale (Scale) and drift (Drift) coefficients:
\begin{align*}
&\;{\bf Scale}\;{\bf 1:} c_{1}(x,\gamma_{1})=\exp\left\{\frac{\gamma_{1,1}+\gamma_{1,2}x+x^{2}}{1+x^{2}}\right\};
\;{\bf Scale}\;{\bf 2:} c_{2}(x,\gamma_{2})=\exp\left\{\frac{\gamma_{2,1}+x+\gamma_{2,3}x^{2}}{1+x^{2}}\right\}; \\
&\;{\bf Scale}\;{\bf 3:} c_{3}(x,\gamma_{3})=\exp\left\{\frac{1+\gamma_{3,2}x+\gamma_{3,3}x^{2}}{1+x^{2}}\right\};
\;{\bf Scale}\;{\bf 4:} c_{4}(x,\gamma_{4})=\exp\left\{\frac{1+\gamma_{4,2}x}{1+x^{2}}\right\}; \\
&\;{\bf Scale}\;{\bf 5:} c_{5}(x,\gamma_{5})=\exp\left\{\frac{1+\gamma_{5,3}x^{2}}{1+x^{2}}\right\}; 
\;{\bf Scale}\;{\bf 6:} c_{6}(x,\gamma_{6})=\exp\left\{\frac{\gamma_{6,2}x+x^{2}}{1+x^{2}}\right\}; \\
&\;{\bf Scale}\;{\bf 7:} c_{7}(x,\gamma_{7})=\exp\left\{\frac{x+\gamma_{7,3}x^{2}}{1+x^{2}}\right\},
\end{align*}
and
\begin{align*}
{\bf Drift}\;{\bf 1:}\; a_{1}(x,\alpha_{1})=-\alpha_{1}(x-1);
\;{\bf Drift}\;{\bf 2:}\; a_{2}(x,\alpha_{2})=-\alpha_{2}x-1;
\;{\bf Drift}\;{\bf 3:}\; a_{3}(x,\alpha_{3})=-\alpha_{3}.
\end{align*}
Each candidate model is given by a combination of these diffusion and drift coefficients; for example, in the case of Scale 1 and Drift 1, we consider the statistical model
\begin{align*}
dX_{t}=-\alpha_{1}(X_{t}-1)dt+\exp\left\{\frac{\gamma_{1,1}+\gamma_{1,2}X_{t}+X_{t}^{2}}{1+X_{t}^{2}}\right\}dw_{t}.
\end{align*}
In this example, although the candidate models do not include the true model, the optimal parameter $(\gamma_{m_{1}}^{\star},\alpha_{m_{2}}^{\star})$ and optimal model indices $m_{1}^{\star}$ and $m_{2}^{\star}$ can be obtained by the functions
\begin{align*}
G_{1}^{(m_{1})}(\gamma_{m_{1}})&=-\int_{\mbbr}\left\{\log c_{m_{1}}(x,\gamma_{m_{1}})^{2}+\frac{1}{c_{m_{1}}(x,\gamma_{m_{1}})^{2}}\right\}\pi_{0}(dx), \\
G_{2}^{(m_{2}|m_{1}^{\star})}(\alpha_{m_{2}})&=-\int_{\mbbr}c_{m_{1}^{\star}}(x,\gamma_{m_{1}^{\star}}^{\star})^{-2}\left\{-\frac{x}{2}-a_{m_{2}}(x,\alpha_{m_{2}})\right\}^{2}\pi_{0}(dx),
\end{align*}
where $\pi_{0}(dx)=\frac{1}{\sqrt{2\pi}}\exp(-x^{2}/2)dx$.
The definition of the optimal model, Tables \ref{simu:tab1}, and \ref{simu:tab2} provide that the optimal model consists of Scale 1 and Drift 1.

Table \ref{simu:tab3} summarizes the comparison results of model selection frequency and the mean of $w_{m_{1},m_{2}}$.
The indicator of the optimal model defined by Scale 1 and Drift 1 is given by $w_{1,1}$.
For all cases, the optimal model is frequently selected,  and the value of $w_{1,1}$ is the highest.
Also observed is that the frequencies that the optimal model is selected and the value of $w_{1,1}$ become higher as $n$ increases.
This result demonstrates the model selection consistency of proposed QBIC.

\subsection{Ergodic L\'{e}vy driven SDE model}
The sample data $\mbX_{n}=(X_{t_{j}})_{j=0}^{n}$ with $t_{j}=jh_{n}$ is obtained from
\begin{align*}
dX_{t}=-\frac{1}{2}X_{t}dt+\frac{1}{1+X_{t-}^{2}}dZ_{t},\quad t\in[0,T_{n}],\quad X_{0}=0,
\end{align*}
where $T_{n}=nh_{n}$, the driving noise process is the normal inverse Gaussian L\'{e}vy process satisfying $\mathcal{L}(Z_{t})=NIG(3,0,3t,0)$.
In this example, we consider the following candidate scale (Scale) and drift (Drift) coefficients:
\begin{align*}
&\;{\bf Scale}\;{\bf 1:} c_{1}(x,\gamma_{1})=\gamma_{1};
\;{\bf Scale}\;{\bf 2:} c_{2}(x,\gamma_{2})=\exp\left\{\frac{1}{2}(\gamma_{2,1}\cos x +\gamma_{2,2}\sin x)\right\}; \\
&\;{\bf Scale}\;{\bf 3:} c_{3}(x,\gamma_{3})=\frac{\gamma_{3}}{1+x^{2}};
\;{\bf Scale}\;{\bf 4:} c_{4}(x,\gamma_{4})=\frac{1+\gamma_{4}x^{2}}{1+x^{2}},
\end{align*}
and
\begin{align*}
{\bf Drift}\;{\bf 1:}\; a_{1}(x,\alpha_{1})=-\alpha_{1,1}x-\alpha_{1,2};
\;{\bf Drift}\;{\bf 2:}\; a_{2}(x,\alpha_{2})=-\alpha_{2}x;
\;{\bf Drift}\;{\bf 3:}\; a_{3}(x,\alpha_{3})=-\alpha_{3}.
\end{align*}
Each candidate model is constructed in a similar manner as Section \ref{subsec_sim1}.
Then, the true model consists of Scale 3 and Drift 2 with $(\gamma_{3},\alpha_{2})=(1,-\frac{1}{2})$. 
Note that Scale 1, 2, and Drift 3 are misspecified coefficients.
From Table \ref{simu:tab4}, we can show that the tendencies of model selection frequency and model weight are analogous to Section \ref{subsec_sim1}.


\section{Appendix}

\noindent\textbf{Proof of Theorem \ref{YU:se}}
Since in the correctly specified and semi-misspecified diffusion cases, Theorem \ref{YU:se} can be shown in a similar way as \cite{EguMas18a}, we consider the cases where the rate of convergence for scale estimator is $\sqrt{T_{n}}$.
For the simplification of the following discussion, we hereafter deal with the zero-extended version of $\mbbg_{1,n}(\gam)$ and $\pi_1(\gam)$: they vanish outside of $\Theta_\gam$.
Applying the change of variable, we have 
\begin{align*}
&\log\left(\int_{\Theta_\gam}\exp\left(\mbbg_{1,n}(\gam)\right)\pi_1(\gam)d\gam\right)\\
&=\mbbg_{1,n}(\ges)-\frac{p_\gam}{2}\log n+\log\left(\int_{\mbbr^{p_\gam}}\exp\left\{\mbbg_{1,n}\left(\ges+\frac{t}{\sqrt{n}}\right)-\mbbg_{1,n}(\ges)\right\}\pi_1\left(\ges+\frac{t}{\sqrt{n}}\right)dt\right).
\end{align*}
Below we show that 
\begin{align}
&\nn\log\left(\int_{\mbbr^{p_\gam}}\exp\left\{\mbbg_{1,n}\left(\ges+\frac{t}{\sqrt{n}}\right)-\mbbg_{1,n}(\ges)\right\}\pi_1\left(\ges+\frac{t}{\sqrt{n}}\right)dt\right)\\
&=\log\pi_1\left(\gam^\star\right)+\frac{p_\gam}{2}\log 2\pi-\frac{1}{2}\log \det \mci_\gam+o_p(1).\label{yu: convp}
\end{align}
It follows from the continuous mapping theorem, \cite[THEOREM 2 (d)]{Fer96} and the estimates of the GQL and GQMLE given in the papers  \cite{Yos92}, \cite{Kes97}, \cite{UchYos11}, \cite{UchYos12}, \cite{Mas13-1}, and \cite{Ueh18} that 
\begin{equation}
\frac{1}{n}\sup_{\gam\in\Theta_\gam}|\p_\gam^3 \mbbg_{1,n}(\gam)|=O_p(1), \label{yu: sb}
\end{equation}
and that for any subsequence $\{n_j\}\subset \{n\}$, we can pick a subsubsequence $\{n_{k_j}\}\subset\{n_j\}$ fulfilling that for all $t\in\mbbr^{p_\gam}$,
\begin{align}
&\left|\frac{1}{n_{k_j}}\p^2_\gam \mbbg_{1,n_{k_j}}\left(\hat{\gam}_{n_{k_j}}\right)+\mci_\gam\right|\asc0,\label{yu:conve1} \\
& \sup_{\gam\in\Theta_\gam}\left|\frac{1}{n_{k_j}}\mbbg_{1,n_{k_j}}(\gam)-\mbbg_1(\gam)\right|\asc0, \label{yu:conve2} \\
&\left|\pi_1\left(\hat{\gam}_{n_{k_j}}+\frac{t}{\sqrt{n_{k_j}}}\right)-\pi_1\left(\gam^\star\right)\right|\asc0, \label{yu:conve3} \\
&\left|\mbbg_1\left(\hat{\gam}_{n_{k_j}}+\frac{t}{\sqrt{n_{k_j}}}\right)-\mbbg_1\left(\gam^\star+\frac{t}{\sqrt{n_{k_j}}}\right)\right|+\left|\mbbg_1(\hat{\gam}_{n_{k_j}})-\mbbg_1(\gam^\star)\right|\asc0. \label{yu:conve4}
\end{align}
For any $\ep>0$, \eqref{yu: sb} enables us to find large enough $N\in\mbbn$ and $M>0$ such that for all $n\geq N$,
\begin{equation*}
P\left(\frac{1}{n}\sup_{\gam\in\Theta_\gam}|\p_\gam^3 \mbbg_{1,n}(\gam)|>M\right)<\ep.
\end{equation*}
Since our main focus here is the convergence in probability \eqref{yu: convp}, we can hereafter suppose that
\begin{equation*}
\frac{1}{n}\sup_{\gam\in\Theta_\gam}|\p_\gam^3 \mbbg_{1,n}(\gam)|\leq M
\end{equation*}
for sufficiently large $n$ and $M$ (and we deal with such $n$ below). 
We also write the set $E\subset\Omega$ on which \eqref{yu:conve1}-\eqref{yu:conve4} hold.
For a positive constant $\del$, we divide $\mbbr^{p_\gam}$ into 
\begin{align*}
&D_{1,n}:=\left\{t\in\mbbr^{p_\gam}: |t|< \del n^{\frac{1}{2}}\right\}, \\
&D_{2,n}:=\left\{t\in\mbbr^{p_\gam}: |t|\geq \del n^{\frac{1}{2}}\right\}.
\end{align*}
For any set $A$, we define the indicator function $\mbbi_{A}(\cdot)$ by:
\begin{align*}
\mbbi_{A}(t)=\begin{cases}1,&t\in A,\\0,& \text{otherwise}.\end{cases}
\end{align*}
First we look at the integration on $D_{1,n}$.
For $t\in\mbbr^{p_\gam}$, we write
\begin{equation*}
R_{n_{k_j}}=\left|\frac{1}{2n_{k_j}}\p^2_\gam \mbbg_{1,n_{k_j}}\left(\hat{\gam}_{n_{k_j}}\right)+\frac{1}{2}\mci_\gam\right|+\del\left|\frac{1}{6n_{k_j}}\int_0^1 \p^3_\gam \mbbg_{1,n_{k_j}}\left(\hat{\gam}_{n_{k_j}}+\frac{t}{\sqrt{n_{k_j}}}u\right)du\right|.
\end{equation*}
For any $\omega\in E$ and $\ep>0$, it follows from \eqref{yu:conve1} and the boundedness of $\frac{1}{n}\sup_{\gam\in\Theta_\gam}|\p_\gam^3 \mbbg_{1,n}(\gam)|$ that there exist $N(\omega)\in\mbbn$ and small enough $\zeta(\omega)>0$ being independent of $t$ satisfying that for all $n_{k_j}\geq N(\omega)$ and $\del<\zeta(\omega)$, 
\begin{equation*}
R_{n_{k_j}}(\omega)\leq \ep.
\end{equation*}
Hence, for all $t\in\mbbr^{p_\gam}$ and $\omega\in E$, by taking $\zeta(\omega)$ and $N(\omega)$ sufficiently small, Taylor's expansion around $\ges$, and \eqref{yu:conve1} yield that for all $\del<\zeta(\omega)$ and $n_{k_j}\geq N(\omega)$, 
\begin{align}
&\nn\exp\left\{\mbbg_{1,n_{k_j}}\left(\hat{\gam}_{n_{k_j}}+\frac{t}{\sqrt{n_{k_j}}}\right)-\mbbg_{1,n_{k_j}}\left(\hat{\gam}_{n_{k_j}}\right)\right\}\pi_1\left(\hat{\gam}_{n_{k_j}}+\frac{t}{\sqrt{n_{k_j}}}\right)\mbbi_{D_{1,n}}(t)\\
&\nn\leq \sup_{\gam\in\Theta_\gam} \pi_1(\gam) \exp\left(-\frac{1}{2}\mci_\gam[t,t]+|t|^2R_{n_{k_j}}\right)\\
&\leq \sup_{\gam\in\Theta_\gam} \pi_1(\gam)  \exp\left\{-\frac{1}{4}\mci_\gam[t,t]\right\} \label{yu:dconve},
\end{align}
and the right-hand-side is integrable over $\mbbr^{p_\gam}$.
The above estimates and \eqref{yu:conve3} also imply that for all $\ep'>0$, $t\in\mbbr^{p_\gam}$, and $\omega\in E$, we can choose $N(\omega)\in\mbbn$ and $\zeta(\omega)>0$ such that for all $n_{k_j}\geq N(\omega)$ and $\del<\zeta(\omega)$, 
\begin{align}
&\nn\left|\exp\left\{\mbbg_{1,n_{k_j}}\left(\hat{\gam}_{n_{k_j}}+\frac{t}{\sqrt{n_{k_j}}}\right)-\mbbg_{1,n_{k_j}}\left(\hat{\gam}_{n_{k_j}}\right)\right\}\pi_1\left(\hat{\gam}_{n_{k_j}}+\frac{t}{\sqrt{n_{k_j}}}\right)-\exp\left(-\frac{1}{2}\mci_\gam [t,t]\right)\pi_1\left(\gam^\star\right)\right|\mbbi_{D_{1,{n_{k_j}}}}(t)\\
&\nn\leq \exp\left(-\frac{1}{2}\mci_\gam [t,t]\right)\left(\left|\pi_1\left(\hat{\gam}_{n_{k_j}}+\frac{t}{\sqrt{n_{k_j}}}\right)-\pi_1\left(\gam^\star\right)\right|+\sup_{\gam\in\Theta_\gam}\pi_1\left(\gam\right)\left|\exp\left(|t|^2R_{n_{k_j}}\right)-1\right| \mbbi_{D_{1,{n_{k_j}}}}(t)\right)\\
&\leq \ep'. \label{yu:pconve}
\end{align}
From \eqref{yu:dconve}, \eqref{yu:pconve} and the dominated convergence theorem, we finally get
\begin{align}
&\nn\log\left(\int_{\mbbr^{p_\gam}} \exp\left\{\mbbg_{1,n_{k_j}}\left(\hat{\gam}_{n_{k_j}}+\frac{t}{\sqrt{n_{k_j}}}\right)-\mbbg_{1,n_{k_j}}\left(\hat{\gam}_{n_{k_j}}\right)\right\}\pi_1\left(\hat{\gam}_{n_{k_j}}+\frac{t}{\sqrt{n_{k_j}}}\right)\mbbi_{D_{1,n_{k_j}}}(t)dt\right)\\
&=\nn\log\Bigg(\int_{\mbbr^{p_\gam}}\Bigg[\exp\left\{\mbbg_{1,n_{k_j}}\left(\hat{\gam}_{n_{k_j}}+\frac{t}{\sqrt{n_{k_j}}}\right)-\mbbg_{1,n_{k_j}}\left(\hat{\gam}_{n_{k_j}}\right)\right\}\pi_1\left(\hat{\gam}_{n_{k_j}}+\frac{t}{\sqrt{n_{k_j}}}\right)\\
&\nn\qquad\qquad\qquad\quad-\exp\left(-\frac{1}{2}\mci_\gam [t,t]\Bigg)\pi_1\left(\gam^\star\right)\Bigg]\mbbi_{D_{1,n_{k_j}}}(t)dt+\int_{\mbbr^{p_\gam}}\exp\left(-\frac{1}{2}\mci_\gam [t,t]\right)\pi_1\left(\gam^\star\right)\mbbi_{D_{1,n_{k_j}}}(t)dt\right)\\
&\nn\asc \log\left(\int_{\mbbr^{p_\gam}}\exp\left(-\frac{1}{2}\mci_\gam [t,t]\right)\pi_1\left(\gam^\star\right)dt\right)\\
&=\log\pi_1\left(\gam^\star\right)+\frac{p_\gam}{2}\log 2\pi-\frac{1}{2}\log \det \mci_\gam.\label{yu:d1}
\end{align}

Now we move on to the evaluation on $D_{2,n}$.
Since
\begin{align*}
&\frac{1}{n_{k_j}}\left(\mbbg_{1,n_{k_j}}\left(\hat{\gam}_{n_{k_j}}+\frac{t}{\sqrt{n_{k_j}}}\right)-\mbbg_{1,n_{k_j}}(\hat{\gam}_{n_{k_j}})\right)\\
&\leq 2\sup_{\gam\in\Theta_\gam}\left|\frac{1}{n_{k_j}}\mbbg_{1,n_{k_j}}(\gam)-\mbbg_1(\gam)\right|+\mbbg_1\left(\gam^\star+\frac{t}{\sqrt{n_{k_j}}}\right)-\mbbg_1(\gam^\star)\\
&+\mbbg_1\left(\hat{\gam}_{n_{k_j}}+\frac{t}{\sqrt{n_{k_j}}}\right)-\mbbg_1\left(\gam^\star+\frac{t}{\sqrt{n_{k_j}}}\right)-\mbbg_1(\hat{\gam}_{n_{k_j}})+\mbbg_1(\gam^\star),
\end{align*}
the identifiability condition, \eqref{yu:conve2}, and \eqref{yu:conve4} imply that on $D_{2,n}$ and for all large enough $n_{k_j}$, there exists a positive constant $\ep''$ satisfying
\begin{equation*}\label{ep}
\frac{1}{n_{k_j}}\left(\mbbg_{1,n_{k_j}}\left(\hat{\gam}_{n_{k_j}}+\frac{t}{\sqrt{n_{k_j}}}\right)-\mbbg_{1,n_{k_j}}\left(\hat{\gam}_{n_{k_j}}\right)\right)<-\ep'',
\end{equation*}
almost surely.
Thus we arrive at 
\begin{align}
&\nn\int_{\mbbr^{p_\gam}} \exp\left\{\mbbg_{1,n_{k_j}}\left(\hat{\gam}_{n_{k_j}}+\frac{t}{\sqrt{n_{k_j}}}\right)-\mbbg_{1,n_{k_j}}\left(\hat{\gam}_{n_{k_j}}\right)\right\}\pi_1\left(\hat{\gam}_{n_{k_j}}+\frac{t}{\sqrt{n_{k_j}}}\right)\mbbi_{D_{2,n_{k_j}}}(t)dt\\
&\leq \exp\left(-n_{k_j}\ep''\right)\int_{\mbbr^{p_\gam}}\pi_1\left(\hat{\gam}_{n_{k_j}}+\frac{t}{\sqrt{n_{k_j}}}\right)dt \asc0.\label{yu:d2}
\end{align}
Again applying the converse of \cite[Theorem 2(d)]{Fer96}, \eqref{yu:d1} and \eqref{yu:d2} imply the desired result.
As for $\log\left(\int_{\Theta_\al}\exp\left(\mbbg_{2,n}(\al)\right)\pi_2(\al)d\al\right)$, the proof is similar, thus we omit its details.

 \medskip
\noindent\textbf{Proof of Theorem \ref{se:thm.modcon}}
For proof of Theorem \ref{se:thm.modcon}, we consider the nested model selection case and non-nested model selection case.
The (non-)nested model means that the candidate models (do not) include the optimal model.
In a similar way as \cite[Theorems 3.3]{EguMas18b} and \cite[Theorems 5.5]{EguMas18a}, we can prove that Theorem \ref{se:thm.modcon} is established for the nested model.
Below, we will deal with the non-nested model selection case.

In the non-nested model selection case, because of the definition of the optimal model, we have $G_{1}^{(m_{1}^{\star})}(\gamma_{m_{1}^{\star}}^{\star})>G_{1}^{(m_{1})}(\gamma_{m_{1}}^{\star})$ a.s. for every $m_{1}\neq m_{1}^{\star}$.
Further, assumptions give the equations
\begin{align*}
\frac{1}{n}G_{1,n}^{(m_{1})}(\hat{\gamma}_{m_{1},n})&=\frac{1}{n}G_{1,n}^{(m_{1})}(\gamma_{m_{1}}^{\star})+o_{p}(1)=G_{1}^{(m_{1})}(\gamma_{m_{1}}^{\star})+o_{p}(1), \\
\frac{1}{n}G_{1,n}^{(m_{1}^{\star})}(\hat{\gamma}_{m_{1}^{\star},n})&=\frac{1}{n}G_{1,n}^{(m_{1}^{\star})}(\gamma_{m_{1}^{\star}}^{\star})+o_{p}(1)=G_{1}^{(m_{1}^{\star})}(\gamma_{m_{1}^{\star}}^{\star})+o_{p}(1).
\end{align*} 
Hence, for any $m_{1}\in\{1,\ldots,M_{1}\}\backslash\{m_{1}^{\star}\}$,
\begin{align}
\mathbb{P}\left(\mathrm{QBIC}_{1,n}^{(m_{1}^{\star})}-\mathrm{QBIC}_{1,n}^{(m_{1})}>0\right)
&=\mathbb{P}\left\{\frac{1}{n}\left(G_{1,n}^{(m_{1}^{\star})}(\hat{\gamma}_{m_{1}^{\star},n})-G_{1,n}^{(m_{1})}(\hat{\gamma}_{m_{1},n})\right)>\left(p_{\gamma_{m_{1}^{\star}}}-p_{\gamma_{m_{1}}}\right)\frac{\log n}{n}\right\} \nn\\
&=\mathbb{P}\left\{G_{1}^{(m_{1}^{\star})}(\gamma_{m_{1}^{\star}}^{\star})-G_{1}^{(m_{1})}(\gamma_{m_{1}}^{\star})>o_{p}(1)\right\} \nn\\
&=\mathbb{P}\left\{G_{1}^{(m_{1}^{\star})}(\gamma_{m_{1}^{\star}}^{\star})-G_{1}^{(m_{1})}(\gamma_{m_{1}}^{\star})>0\right\}+o(1) \nn\\
&\to1 \label{se:prf.thm.modcon1}
\end{align}
as $n\to\infty$.
As with \eqref{se:prf.thm.modcon1}, we can show that for any $m_{2}\in\{1,\ldots,M_{2}\}\backslash\{m_{2}^{\star}\}$,
\begin{align}
\mathbb{P}\left(\mathrm{QBIC}_{2,n}^{(m_{2}^{\star}|m_{1}^{\star})}-\mathrm{QBIC}_{2,n}^{(m_{2}|m_{1}^{\star})}>0\right)
&\to1. \label{se:prf.thm.modcon2}
\end{align}
From \eqref{se:prf.thm.modcon1} and \eqref{se:prf.thm.modcon2}, we have
\begin{align*}
\mathbb{P}\left(\mathrm{QBIC}_{2,n}^{(m_{2}^{\star}|\hat{m}_{1,n})}-\mathrm{QBIC}_{2,n}^{(m_{2}|\hat{m}_{1,n})}>0\right)&=\mathbb{P}\left(\mathrm{QBIC}_{2,n}^{(m_{2}^{\star}|\hat{m}_{1,n})}-\mathrm{QBIC}_{2,n}^{(m_{2}|\hat{m}_{1,n})}>0, \hat{m}_{1,n}=m_{1}^{\star}\right) \\
&\quad+\mathbb{P}\left(\mathrm{QBIC}_{2,n}^{(m_{2}^{\star}|\hat{m}_{1,n})}-\mathrm{QBIC}_{2,n}^{(m_{2}|\hat{m}_{1,n})}>0, \hat{m}_{1,n}\neq m_{1}^{\star}\right) \\
&\leq\mathbb{P}\left(\mathrm{QBIC}_{2,n}^{(m_{2}^{\star}|m_{1}^{\star})}-\mathrm{QBIC}_{2,n}^{(m_{2}|m_{1}^{\star})}>0\right) \\
&\quad+\mathbb{P}\left(\hat{m}_{1,n}\neq m_{1}^{\star}\right) \\
&=\mathbb{P}\left(\mathrm{QBIC}_{2,n}^{(m_{2}^{\star}|m_{1}^{\star})}-\mathrm{QBIC}_{2,n}^{(m_{2}|m_{1}^{\star})}>0\right) \\
&\quad+\mathbb{P}\left(\mathrm{QBIC}_{1,n}^{(m_{1}^{\star})}<\mathrm{QBIC}_{1,n}^{(m_{1}^{\star})}\right) \\
&\to1+0=1.
\end{align*}
The proof of Theorem \ref{se:thm.modcon} is complete.

\subsection*{Acknowledgement}
The author wishes to thank the associate editor and the two anonymous referees for careful reading and valuable comments which helped to greatly improve the paper.
This work was supported by JST CREST Grant Number JPMJCR14D7, Japan.

\clearpage

\begin{table}[t]
\begin{center}
\caption{The values of $G_{1}^{(m_{1})}(\gamma_{m_{1}}^{\star})$ for each candidate diffusion coefficient.}
\begin{tabular}{r r r r r r r r r} \hline
& & Scale 1 & Scale 2 & Scale 3 & Scale 4 & Scale 5 & Scale 6 & Scale 7 \\ \hline
& & & & & & & & \\[-3mm]
$G_{1}^{(m_{1})}(\gamma_{m_{1}}^{\star})$ & & -1.2089 & -1.2822 & -1.4833 & -1.6225 & -1.4833 & -1.2602 & -3.2860 \\[1mm] \hline
\end{tabular}
\label{simu:tab1}
\end{center}
\end{table}

\begin{table}[t]
\begin{center}
\caption{The values of $G_{2}^{(m_{2}|m_{1}^{\star})}(\alpha_{m_{2}}^{\star})$ for each candidate drift coefficient.}
\begin{tabular}{r r r r r} \hline
& & Drift 1 & Drift 2 & Drift 3 \\ \hline
& & & & \\[-3mm]
$G_{2}^{(m_{2}|m_{1}^{\star})}(\alpha_{m_{2}}^{\star})$ & & -0.0624 & -0.8193 & -0.0979 \\[1mm] \hline
\end{tabular}
\label{simu:tab2}
\end{center}
\end{table}

\clearpage

\begin{table}[t]
\begin{center}
\caption{The mean of model weight $w_{m_{1},m_{2}}$ and model selection frequencies for various situations. The optimal model consists of Scale 1 and Drift 1.}
\begin{tabular}{r r l l r r r r r r r} \hline
$T_{n}$ & $h_{n}$ & & & Scale $1^{\ast}$ & Scale 2 & Scale 3 & Scale 4 & Scale 5 & Scale 6 & Scale 7 \\ \hline
10 & 0.01 & Drift $1^{\ast}$ & frequency & \bf{409} & 72 & 5 & 1 & 5 & 95 & 70 \\
\multicolumn{2}{r}{$(n=1000)$} & & weight & \bf{30.27} & 7.26 & 0.41 & 0.04 & 0.41 & 7.57 & 5.38 \\
& & Drift 2 & frequency & 60 & 84 & 2 & 0 & 0 & 31 & 22 \\
& & & weight & 5.94 & 6.67 & 0.13 & 0.01 & 0.02 & 2.64 & 1.98 \\
& & Drift 3 & frequency & 125 & 5 & 0 & 0 & 0 & 3 & 11 \\
& & & weight & 22.91 & 2.50 & 0.15 & 0.04 & 0.10 & 3.06 & 2.51 \\ \hline
10& 0.005 & Drift $1^{\ast}$ & frequency & \bf{449} & 86 & 6 & 0 & 4 & 73 & 45 \\
\multicolumn{2}{r}{$(n=2000)$} & & weight & \bf{33.19} & 8.07 & 0.53 & 0.02 & 0.30 & 5.61 & 3.51 \\
& & Drift 2 & frequency & 64 & 96 & 3 & 0 & 0 & 26 & 8 \\
& & & weight & 6.61 & 7.65 & 0.19 & 0.00 & 0.01 & 1.95 & 0.89 \\
& & Drift 3 & frequency & 129 & 4 & 1 & 0 & 0 & 2 & 4 \\ 
& & & weight & 24.63 & 2.94 & 0.26 & 0.02 & 0.07 & 2.07 & 1.48 \\ \hline
50 & 0.01 &Drift $1^{\ast}$ & frequency & \bf{832} & 58 & 2 & 0 & 1 & 1 & 12 \\
\multicolumn{2}{r}{$(n=5000)$} & & weight & \bf{62.59} & 5.19 & 0.19 & 0.00 & 0.10 & 0.08 & 0.99 \\
& & Drift 2 & frequency & 2 & 13 & 0 & 0 & 0 & 0 & 0 \\
& & & weight & 0.29 & 1.12 & 0.00 & 0.00 & 0.00 & 0.00 & 0.04 \\
& & Drift 3 & frequency & 79 & 0 & 0 & 0 & 0 & 0 & 0 \\
& & & weight & 28.43 & 0.74 & 0.01 & 0.00 & 0.00 & 0.01 & 0.21 \\ \hline
50 & 0.005 & Drift $1^{\ast}$ & frequency & \bf{841} & 59 & 3 & 0 & 2 & 0 & 7 \\
\multicolumn{2}{r}{$(n=10000)$} & & weight & \bf{62.80} & 5.30 & 0.30 & 0.00 & 0.19 & 0.00 & 0.59 \\
& & Drift 2 & frequency & 3 & 13 & 0 & 0 & 0 & 0 & 0 \\
& & & weight & 0.31 & 1.15 & 0.00 & 0.00 & 0.00 & 0.00 & 0.00 \\
& & Drift 3 & frequency & 72 & 0 & 0 & 0 & 0 & 0 & 0 \\
& & & weight & 28.46 & 0.76 & 0.01 & 0.00 & 0.00 & 0.00 & 0.12 \\ \hline
\end{tabular}
\label{simu:tab3}
\end{center}
\end{table}

\clearpage

\begin{table}[t]
\begin{center}
\caption{The mean of model weight $w_{m_{1},m_{2}}$ and model selection frequencies for various situations. The true model consists of Scale 3 and Drift 2.}
\begin{tabular}{r r l l r r r r} \hline
$T_{n}$ & $h_{n}$ & & & Scale 1 & Scale 2 & Scale $3^{\ast}$ & Scale 4 \\ \hline
10 & 0.01 & Drift 1 & frequency & 3 & 38 & 169 & 27 \\
\multicolumn{2}{r}{$(n=1000)$} & & weight & 0.54 & 5.07 & 25.09 & 7.94 \\
& & Drift $2^{\ast}$ & frequency & 12 & 72 & {\bf 548} & 131 \\
& & & weight & 0.91 & 5.55 & {\bf 39.94} & 13.28 \\
& & Drift 3 & frequency & 0 & 0 & 0 & 0 \\
& & & weight & 0.10 & 0.50 & 0.78 & 0.30 \\ \hline
10& 0.005 & Drift 1 & frequency & 1 & 36 & 174 & 28 \\
\multicolumn{2}{r}{$(n=2000)$} & & weight & 0.38 & 4.84 & 26.64 & 6.94 \\
& & Drift $2^{\ast}$ & frequency & 11 & 70 & {\bf 557} & 123 \\
& & & weight & 0.81 & 5.29 & {\bf 42.01} & 11.51 \\
& & Drift 3 & frequency & 0 & 0 & 0 & 0 \\
& & & weight & 0.07 & 0.45 & 0.82 & 0.24 \\ \hline
50 & 0.01 &Drift 1 & frequency & 0 & 0 & 68 & 24 \\
\multicolumn{2}{r}{$(n=5000)$} & & weight & 0.00 & 0.01 & 14.44 & 6.70 \\
& & Drift $2^{\ast}$ & frequency & 0 & 1 & {\bf 659} & 248 \\
& & & weight & 0.00 & 0.09 & {\bf 54.09} & 24.66 \\
& & Drift 3 & frequency & 0 & 0 & 0 & 0 \\
& & & weight & 0.00 & 0.00 & 0.00 & 0.00 \\ \hline
50 & 0.005 & Drift 1 & frequency & 0 & 0 & 69 & 20 \\
\multicolumn{2}{r}{$(n=10000)$} & & weight & 0.00 & 0.01 & 15.31 & 5.84 \\
& & Drift $2^{\ast}$ & frequency & 0 & 1 & {\bf 684} & 226\\
& & & weight & 0.00 & 0.09 & {\bf 57.41} & 21.35 \\
& & Drift 3 & frequency & 0 & 0 & 0 & 0 \\
& & & weight & 0.00 & 0.00 & 0.00 & 0.00 \\\hline
\end{tabular}
\label{simu:tab4}
\end{center}
\end{table}

\clearpage

\bibliographystyle{abbrv}

\end{document}